\documentclass[10pt]{article}
 \font\smallit=cmti10

\usepackage{amssymb,latexsym,amsmath,epsfig,amsthm,color,bm} %% Add other packages as necessary

\makeatletter

\renewcommand{\@seccntformat}[1]{\csname the#1\endcsname. }

\makeatother

 \newtheorem{theorem}{Theorem}[section]
 \newtheorem{lemma}[theorem]{Lemma}
 
 \newtheorem{proposition}[theorem]{Proposition}
 \newtheorem{corollary}[theorem]{Corollary}
 \newtheorem{definition}[theorem]{Definition}
 \newtheorem{remark}[theorem]{Remark}
 \newtheorem{example}[theorem]{Example}
\begin{document}
\begin{center}
 {\bf An Explicit Construction of Orthogonal Basis in $p$-adic Fields}
 \vskip 30pt

 {\bf Chi Zhang and Yingpu Deng}\\
 {\smallit State Key Laboratory of Mathematical Sciences, Academy of Mathematics and Systems Science, Chinese Academy of Sciences, Beijing 100190, People's Republic of China}\\
 {and}\\
 {\smallit School of Mathematical Sciences, University of Chinese Academy of Sciences, Beijing 100049, People's Republic of China}\\

 \vskip 10pt

 {\tt zhangchi171@mails.ucas.ac.cn, dengyp@amss.ac.cn}\\

 \end{center}
 \vskip 30pt

\centerline{\bf Abstract}
In 2021, the $p$-adic signature scheme and public-key encryption cryptosystem were introduced. These schemes have good efficiency but are shown to be not secure. The attack succeeds because the extension fields used in these schemes are totally ramified. In order to avoid this attack, the extension field should have a large residue degree. In this paper, we propose a method of constructing a kind of specific orthogonal basis in $p$-adic fields with a large residue degree, which would be helpful to modify the $p$-adic signature scheme and public-key encryption cryptosystem.

\vskip 10pt
2010 Mathematics Subject Classification: Primary 11F85, Secondary 94A60.\par
Key words and phrases: Local field, Orthogonal basis, $p$-adic lattice.

\noindent

\pagestyle{myheadings}

 \thispagestyle{empty}
 \baselineskip=12.875pt
 \vskip 20pt

\section{Introduction}

Since Peter Shor \cite{ref-11} proved that the classical public-key cryptosystems such as RSA and ElGamal would be broken by future quantum computer, researchers have been dedicated to finding cryptographic primitives which are quantum-resistant. In 2022, NIST \cite{ref-16} announced four algorithms which passed the third round of post-quantum cryptography standardization solicitation and began the fourth round. They are CRYSTALS-Kyber \cite{ref-12}, CRYSTALS-Dilithium \cite{ref-13}, Falcon \cite{ref-14} and SPHINCS$^{+}$ \cite{ref-15}. Three of them are lattice-based and one of them is hash-based. The lack of diversity among post-quantum assumptions is widely recognized as a big, open issue in the field. Therefore, finding new  post-quantum assumptions is of vital significance.\par
The $p$-adic numbers $\mathbb{Q}_p$ were invented by Hensel in the late 19th century. The concept of a local field is an abstraction of the field $\mathbb{Q}_p$.  Local fields provide a natural tool to solve many number-theoretic problems. They are ubiquitous in modern algebraic number theory and arithmetic geometry. Lattices can also be defined in local fields such as $p$-adic fields, see \cite{ref-5}. Interestingly, $p$-adic lattices possess some properties which lattices in Euclidean spaces do not have, see \cite{ref-6}. However, applications of $p$-adic lattices in cryptography were developed only recently.\par
In 2021, by introducing a trapdoor function with an orthogonal basis of a $p$-adic lattice, Deng et al. \cite{ref-2} constructed the first signature scheme and public-key encryption cryptosystem based on $p$-adic lattices. As the $p$-adic analogues of the lattices in Euclidean spaces, it is reasonable to expect hard problems in $p$-adic lattices to be quantum-resistant, which might provide new alternative candidates to construct post-quantum cryptographic primitives.\par
The experimental results \cite{ref-2} demonstrated that the new schemes achieve good efficiency. As for security, Zhang \cite{ref-7} found that these schemes are not secure because the extension fields used in these schemes are totally ramified. In order to avoid this attack, he suggested that the extension field should have a large residue degree.\par
In a totally ramified extension field $K/\mathbb{Q}_p$, a uniformizer $\pi$ generates an orthogonal basis of $K$. But in a general extension field $K/\mathbb{Q}_p$, we can not find an orthogonal basis of $K$ as easily as in a totally ramified extension field. Therefore, the crucial point of such a scheme is to construct an orthogonal basis of $K$.\par
Given a extension field $K$ over $\mathbb{Q}_p$ of degree $n$, we can use the Round 2 Algorithm \cite{ref-18} or the Round 4 Algorithm \cite{ref-19} to obtain a basis of the maximal order $\mathcal{O}_K$ and then compute its orthogonal basis. However, these algorithms involve computation of large matrices. They require storage of the order of $n^3$ in the worst case.\par
In order to reduce the storage requirement, we consider the problem from another perspective. Instead of trying computing the maximal order, we construct an orthogonal basis directly and then compute the extention field it generates. The storage requirement of this method is of the order of $n^2$ in the worst case.\par
This paper is organized as follows. In Section \ref{se-2}, we recall some basic definitions. In Section \ref{se-3}, we give an equivalent condition for orthogonal basis in the extension field of $\mathbb{Q}_p$. Then, we construct a kind of specific orthogonal basis in Section \ref{se-4} and realize it with roots of unity in Section \ref{se-5}.

\section{Preliminaries}\label{se-2}

In this section, we recall some basic facts about local fields and $p$-adic lattices. More details about local fields can be found in \cite{ref-4,ref-4.5}.\par

\subsection{Norm and Orthogonal Basis}

Let $p$ be a prime. Let $V$ be a vector space over $\mathbb{Q}_p$. A norm $\left\|\cdot\right\|$ on $V$ is a function
$$\left\|\cdot\right\|:V\rightarrow\mathbb{R}$$
such that:

\begin{enumerate}
\item $\left\|{\bm v}\right\|\ge0$ for any ${\bm v}\in V$, and $\left\|{\bm v}\right\|=0$ if and only if ${\bm v}=0$;
\item $\left\|x{\bm v}\right\|=\left|x\right|_{p}\cdot\left\|{\bm v}\right\|$ for any $x\in\mathbb{Q}_p$ and ${\bm v}\in V$;
\item $\left\|{\bm v}+{\bm w}\right\|\le\max{\left\{\left\|{\bm v}\right\|,\left\|{\bm w}\right\|\right\}}$ for any ${\bm v},{\bm w}\in V$.
\end{enumerate}

\noindent Here, $\left|x\right|_{p}$ is the $p$-adic absolute value for a $p$-adic number $x\in\mathbb{Q}_p$.\par
If $\left\|\cdot\right\|$ is a norm on $V$, and if $\left\|{\bm v}\right\|\ne\left\|{\bm w}\right\|$ for ${\bm v},{\bm w}\in V$, then one can prove $\left\|{\bm v}+{\bm w}\right\|=\max{\left\{\left\|{\bm v}\right\|,\left\|{\bm w}\right\|\right\}}$. Weil (\cite{ref-5} page 26) proved the following proposition:

\begin{proposition}[\cite{ref-5}]
Let $V$ be a vector space over $\mathbb{Q}_p$ of finite dimension $n>0$, and let $\left\|\cdot\right\|$ be a norm on $V$. Then there is a decomposition $V=V_1+V_2+\cdots+V_n$ of $V$ into a direct sum of subspaces $V_i$ of dimension $1$, such that
$$\left\|\sum_{i=1}^{n}{{\bm v}_i}\right\|=\max_{1\le i\le n}{\left\|{\bm v}_i\right\|}$$
for any ${\bm v}_i\in V_i$, $i=1,2,\dots,n$.
\end{proposition}

Thus, we can define the orthogonal basis.

\begin{definition}[orthogonal basis]
Let $V$ be a vector space over $\mathbb{Q}_p$ of finite dimension $n>0$, and let $\left\|\cdot\right\|$ be a norm on $V$. We call ${\bm \alpha}_1,{\bm \alpha}_2,\dots,{\bm \alpha}_n$ an orthogonal basis of $V$ over $\mathbb{Q}_p$ if $V$ can be decomposed into the direct sum of $n$ $1$-dimensional subspaces $V_i$'s $(1\le i\le n)$, such that
$$\left\|\sum_{i=1}^{n}{{\bm v}_i}\right\|=\max_{1\le i\le n}{\left\|{\bm v}_i\right\|}$$
for any ${\bm v}_i\in V_i$, $i=1,2,\dots,n$, where $V_i$ is spanned by ${\bm \alpha}_i$. Two subspaces $U$, $W$ of $V$ are said to be orthogonal if the sum $U+W$ is a direct sum and it holds that $\left\|{\bm u}+{\bm w}\right\|=\max\left\{\left\|{\bm u}\right\|, \left\|{\bm w}\right\|\right\}$ for all ${\bm u}\in U$, ${\bm w}\in W$. Actually, the former assumption can be deduced from the latter assumption.
\end{definition}

\subsection{Residue Degree and Ramification Index}

Let $K$ be a finite extension of degree $n$ over the field $\mathbb{Q}_p$ of $p$-adic numbers. Hence $K$ is locally compact and complete. Let us choose an element $\pi$ of maximal absolute value smaller than $1$ and call it a uniformizer. Let
$$R=\left\{x\in K\big|\left|x\right|\le1\right\}$$
and its maximal ideal $P=\pi R$, where $\left|\cdot\right|$ denotes the unique extension of the $p$-adic absolute value to the field $K$ (see \cite{ref-4} page 94). If we consider $K$ as a finite-dimensional vector space over $\mathbb{Q}_p$, then the extended absolute value is also a norm on this vector space. For this absolute-value norm, the second property of norm holds not only for $x\in\mathbb{Q}_p$, bug also for $x\in K$. The residue field $k=R/P$ is finite, hence a finite extension over $\mathbb{F}_p=\mathbb{Z}_p/p\mathbb{Z}_p$.

\begin{definition}[residue degree and ramification index]
The residue degree of the finite extension $K$ over $\mathbb{Q}_p$ is the integer
$$f=\left[k:\mathbb{F}_p\right]=\dim_{\mathbb{F}_p}(k).$$
The ramification index of $K$ over $\mathbb{Q}_p$ is the integer
$$e=\left[\left|K^*\right|:\left|\mathbb{Q}_p^*\right|\right]=\left[\left|K^*\right|:\left|p^{\mathbb{Z}}\right|\right]=\#\left(\left|K^*\right|/p^{\mathbb{Z}}\right),$$
where $\left|K^*\right|$ denotes the value group of $K^*$.
\end{definition}

Similar to the finite extension field over the rational number field $\mathbb{Q}$, we have the following theorem.

\begin{theorem}[\cite{ref-4} page 99]\label{th-0.4}
For each finite extension $K$ over $\mathbb{Q}_p$, we have
$$ef=\left[K:\mathbb{Q}_p\right]=n.$$
\end{theorem}

\subsection{Lattice in $p$-adic Fields}

As in the previous subsection, let $p$ be a prime number, and let $K$ be an extension field of $\mathbb{Q}_p$ of degree $n$. Let $m$ be a positive integer with $1\le m\le n$. Let ${\bm \alpha}_1,{\bm \alpha}_{2},\dots,{\bm \alpha}_m\in K$ be $m$ $\mathbb{Q}_p$-linearly independent vectors. A lattice in $K$ is the set
$$\mathcal{L}({\bm \alpha}_1,{\bm \alpha}_{2},\dots,{\bm \alpha}_m)=\left\{\sum^{m}_{i=1}{a_i{\bm \alpha}_i}\Bigg|a_i\in\mathbb{Z}_p,1\le i\le m\right\}$$
of all $\mathbb{Z}_p$-linear combinations of ${\bm \alpha}_1,{\bm \alpha}_{2},\dots,{\bm \alpha}_m$. The sequence of vectors ${\bm \alpha}_1,{\bm \alpha}_{2},\dots,{\bm \alpha}_m$ is called a basis of the lattice $\mathcal{L}({\bm \alpha}_1,{\bm \alpha}_{2},\dots,{\bm \alpha}_m)$. The integer $m$ is called the rank of the lattice, respectively. When $n=m$, we say that the lattice is of full rank.

\section{An Equivalent Condition for Orthogonal Basis}\label{se-3}

In this section, we give an equivalent condition for orthogonal basis in the extension field of $\mathbb{Q}_p$. We begin with some simple lemmas.

\begin{lemma}\label{le-1.1}
Let $V$ be a vector space over $\mathbb{Q}_p$ of finite dimension $n>0$. Let ${\bm \alpha}_{1},{\bm \alpha}_{2},\dots,{\bm \alpha}_{n}$ be a basis of $V$ over $\mathbb{Q}_p$. Then ${\bm \alpha}_{1},{\bm \alpha}_{2},\dots,{\bm \alpha}_{n}$ is an orthogonal basis of $V$ over $\mathbb{Q}_p$ if and only if 
\begin{equation}\nonumber
\left\|\sum^{n}_{i=1}a_{i}{\bm \alpha}_{i}\right\|=\max_{1\le i\le n}\left\| a_{i}{\bm \alpha}_{i}\right\|
\end{equation}
for all $a_{i}\in\mathbb{Z}_p$, $1\le i\le n$, where at least one $a_{i}\in\mathbb{Z}_p-p\mathbb{Z}_p$.
\begin{proof}
Necessity is immediately from the definition of the orthogonal basis. Now we prove the sufficiency. If there is any $a_{i}=0$, we can just ignore it. So we may assume that $a_{i}\in\mathbb{Q}_p$ and $a_{i}\ne0$, $1\le i\le n$. For $x\in\mathbb{Q}_p$ and $x\ne0$, write $x=p^tu$ with a unit $u\in\mathbb{Z}_p^{*}$, we define ${\rm ord}(x)=t$. Let $s=\min_{1\le i\le n}{\{{\rm ord}(a_{i})\}}$ so that $p^{-s}a_{i}\in\mathbb{Z}_p$ for all $1\le i\le n$ and at least one $a_{i}\in\mathbb{Z}_p-p\mathbb{Z}_p$. Then,
\begin{equation}\nonumber
\left|p^{-s}\right|_{p}\cdot\left\|\sum^{n}_{i=1}a_{i}{\bm \alpha}_{i}\right\|=\left\|\sum^{n}_{i=1}p^{-s}a_{i}{\bm \alpha}_{i}\right\|=\max_{1\le i\le n}\left\| p^{-s}a_{i}{\bm \alpha}_{i}\right\|=\left|p^{-s}\right|_{p}\cdot\max_{1\le i\le n}\left\| a_{i}{\bm \alpha}_{i}\right\|.
\end{equation}
Therefore
\begin{equation}\nonumber
\left\|\sum^{n}_{i=1}a_{i}{\bm \alpha}_{i}\right\|=\max_{1\le i\le n}\left\| a_{i}{\bm \alpha}_{i}\right\|
\end{equation}
for all $a_{i}\in\mathbb{Q}_p$, $1\le i\le n$. Hence ${\bm \alpha}_{1},{\bm \alpha}_{2},\dots,{\bm \alpha}_{n}$ is an orthogonal basis of $V$ over $\mathbb{Q}_p$.
\end{proof}
\end{lemma}

\begin{lemma}\label{le-1.2}
Let $V$ be a vector space over $\mathbb{Q}_p$ of finite dimension $n>0$. Let ${\bm \alpha}_{1},{\bm \alpha}_{2},\dots,{\bm \alpha}_{n}$ be a basis of $V$ over $\mathbb{Q}_p$ such that $\left\|{\bm \alpha}_{1}\right\|=\left\|{\bm \alpha}_{2}\right\|=\cdots=\left\|{\bm \alpha}_{n}\right\|$. Let $\lambda_{1}=\left\|{\bm \alpha}_{1}\right\|$. Then ${\bm \alpha}_{1},{\bm \alpha}_{2},\dots,{\bm \alpha}_{n}$ is an orthogonal basis of $V$ over $\mathbb{Q}_p$ if and only if for all $a_1,a_2,\dots,a_n\in\mathbb{Z}_p$ such that
$$\left\|\sum^{n}_{i=1}a_{i}{\bm \alpha}_{i}\right\|<\lambda_{1},$$
we have $p|a_{i}$ for $1\le i\le n$.
\begin{proof}
Assume ${\bm \alpha}_{1},{\bm \alpha}_{2},\dots,{\bm \alpha}_{n}$ is an orthogonal basis of $V$ over $\mathbb{Q}_p$. Then 
\begin{equation}\nonumber
\left\|\sum^{n}_{i=1}a_{i}{\bm \alpha}_{i}\right\|=\max_{1\le i\le n}\left\| a_{i}{\bm \alpha}_{i}\right\|=\lambda_{1}\cdot\max_{1\le i\le n}\left\| a_{i}\right\|.
\end{equation}
Therefore $\left\|\sum^{n}_{i=1}a_{i}{\bm \alpha}_{i}\right\|<\lambda_{1}$ implies that $\max_{1\le i\le n}\left|a_{i}\right|_{p}<1$. Hence $p|a_{i}$ for all $1\le i\le n$.\par
Conversely, if ${\bm \alpha}_{1},{\bm \alpha}_{2},\dots,{\bm \alpha}_{n}$ is not an orthogonal basis of $V$ over $\mathbb{Q}_p$, then by Lemma \ref{le-1.1}, there exists $a_{i}\in\mathbb{Z}_p$, $1\le i\le n$, and at least one $a_{i}\in\mathbb{Z}_p-p\mathbb{Z}_p$, such that 
\begin{equation}\nonumber
\left\|\sum^{n}_{i=1}a_{i}{\bm \alpha}_{i}\right\|<\max_{1\le i\le n}\left\| a_{i}{\bm \alpha}_{i}\right\|=\lambda_{1},
\end{equation}
which contradicts the assumption of sufficiency in this lemma.
\end{proof}
\end{lemma}

Then we can prove our main theorem of this section. From now on, we use $\left|\cdot\right|$ to denote the $p$-adic absolute value on an extension field $K$ over $\mathbb{Q}_p$ for simplicity. The following theorem is a well-known result (see \cite{ref-10} page 167, Exercise 5A), showing the relation between orthogonality and linear independence. For the sake of completeness, we provide a proof here.

\begin{theorem}\label{th-1.3}
Let $K$ be an extension field of degree $n$ over $\mathbb{Q}_p$. Let $V$ be a subspace of $K$ over $\mathbb{Q}_p$. Assume that ${\bm \alpha}_{1},{\bm \alpha}_{2},\dots,{\bm \alpha}_{m}$ $(m\le n)$ is a basis of $V$ over $\mathbb{Q}_p$ and $\left|{\bm \alpha}_{1}\right|=\left|{\bm \alpha}_{2}\right|=\cdots=\left|{\bm \alpha}_{m}\right|$. Let $\lambda_{1}=\left|{\bm \alpha}_{1}\right|$. Let $\pi$ be a uniformizer of $K$, so there is an integer $s$ such that $\left|\pi^{s}\right|=\lambda_{1}$. Then ${\bm \alpha}_{1},{\bm \alpha}_{2},\dots,{\bm \alpha}_{m}$ is an orthogonal basis of $V$ over $\mathbb{Q}_p$ if and only if $\overline{{\bm \alpha}_{1}},\overline{{\bm \alpha}_{2}},\dots,\overline{{\bm \alpha}_{m}}$ are linearly independent over $\mathbb{F}_p$, where $\overline{{\bm \alpha}_{i}}$ is the image of $\pi^{-s}\cdot{\bm \alpha}_{i}$ in $k=R/P$.
\begin{proof}
By Lemma \ref{le-1.2}, ${\bm \alpha}_{1},{\bm \alpha}_{2},\dots,{\bm \alpha}_{m}$ is an orthogonal basis of $V$ over $\mathbb{Q}_p$ if and only if
$$\left|\sum^{m}_{i=1}a_{i}{\bm \alpha}_{i}\right|<\lambda_{1},\ a_i\in\mathbb{Z}_p \ \Rightarrow\  p|a_{i} \mbox{ for all } 1\le i\le m.$$
Recall that $\left|ab\right|=\left|a\right|\left|b\right|$ for $a,b\in K$. The above statement is equivalent to
$$\left|\sum^{m}_{i=1}a_{i}\pi^{-s}{\bm \alpha}_{i}\right|<1,\ a_i\in\mathbb{Z}_p \ \Rightarrow\ p|a_{i} \mbox{ for all } 1\le i\le m.$$
And this is equivalent to
$$\sum^{m}_{i=1}\overline{a_{i}{\bm \alpha}_{i}}=\overline0,\ a_i\in\mathbb{Z}_p \ \Rightarrow\ \overline{a_{i}}=\overline0 \mbox{ for all } 1\le i\le m,$$
which is the definition of linear independence of $\overline{{\bm \alpha}_{1}},\overline{{\bm \alpha}_{2}},\dots,\overline{{\bm \alpha}_{m}}$ over $\mathbb{F}_p$.
\end{proof}
\end{theorem}

\begin{remark}
In general, it is not efficient to determine orthogonal basis by this theorem. However, in some special cases, it will be useful. For example, in a totally ramified extension field, since $k=\mathbb{F}_p$, any two vectors ${\bm \alpha}_1$ and ${\bm \alpha}_2$ such that $\left|{\bm \alpha}_{1}\right|=\left|{\bm \alpha}_{2}\right|$ can not be extended to an orthogonal basis of $V$ over $\mathbb{Q}_p$.
\end{remark}

\section{Construction of Specific Orthogonal Basis}\label{se-4}

In this section, we use Theorem \ref{th-1.3} to construct a kind of specific orthogonal basis.

\begin{lemma}\label{le-2.1}
Let $K=\mathbb{Q}_p(\theta)$ be an extension field of degree $n$ over $\mathbb{Q}_p$ with $\left|\theta\right|=1$. Let $F$ be the minimal polynomial of $\theta$ over $\mathbb{Q}_p$. Assume that $F$ is reducible modulo $p$. Then $1,\theta,\theta^2,\dots,\theta^{n-1}$ is not an orthogonal basis of $K$ over $\mathbb{Q}_p$.
\begin{proof}
Let $F\equiv gh\pmod{p}$ where $g$ and $h$ have coefficients in $\mathbb{F}_p$. Since $g(\theta)h(\theta)\in pR\subset P$, at least one of $g(\theta)$ and $h(\theta)$ is in $P$. We may assume that $g(\theta)\in P$. Then
 $\overline{g(\theta)}=\overline0$. Since $\deg(g)\le n-1$, $\overline1,\overline\theta,\overline{\theta^2},\dots,\overline{\theta^{n-1}}$ are linearly dependent over $\mathbb{F}_p$. By Theorem \ref{th-1.3}, $1,\theta,\theta^2,\dots,\theta^{n-1}$ is not an orthogonal basis of $K$ over $\mathbb{Q}_p$.
\end{proof}
\end{lemma}

\begin{example}
Let $\theta$ be a primitive $p^l$th root of unity. Then $K=\mathbb{Q}_p(\theta)$ is a totally ramified extension field of degree $n=\varphi(p^l)=p^{l-1}(p-1)$. Since $X^{p^l}-1\equiv (X-1)^{p^l}\pmod{p}$ and the minimal polynomial of $\theta$ is a factor of $X^{p^l}-1$, by Lemma \ref{le-2.1}, $1,\theta,\theta^2,\dots,\theta^{n-1}$ is not an orthogonal basis of $K$ over $\mathbb{Q}_p$. Moreover, we can deduce from the proof of Lemma \ref{le-2.1} that $\left|\theta-1\right|<1$. In fact, $\left|\theta-1\right|=\left|p\right|^{\frac{1}{\varphi(p^l)}}$.
\end{example}

The converse proposition is also true in an unramified extension field, so we have the following theorem.

\begin{theorem}\label{th-2.3}
Let $K=\mathbb{Q}_p(\theta)$ be an unramified extension field of degree $n$ over $\mathbb{Q}_p$ with $\left|\theta\right|=1$. Let $F$ be the minimal polynomial of $\theta$ over $\mathbb{Q}_p$. Then $1,\theta,\theta^2,\dots,\theta^{n-1}$ is an orthogonal basis of $K$ over $\mathbb{Q}_p$ if and only if $F$ is irreducible modulo $p$. 
\begin{proof}
Since $1,\theta,\theta^2,\dots,\theta^{n-1}$ is a basis of $K$ over $\mathbb{Q}_p$, by Theorem \ref{th-1.3}, $1,\theta,\theta^2,\dots,\theta^{n-1}$ is not an orthogonal basis of $K$ over $\mathbb{Q}_p$ if and only if $\overline1,\overline\theta,\overline{\theta^2},\dots,\overline{\theta^{n-1}}$ are linearly dependent over $\mathbb{F}_p$. So there is a polynomial
$$g(X)=\sum^{n-1}_{i=0}a_{i}X^{i}$$
with coefficients $a_{i}\in\mathbb{Z}_p$ for $0\le i\le n-1$, such that $\overline{g(X)}\ne\overline0$ and $\overline{g(\theta)}=\overline0$. Let $G$ be the minimal polynomial of $\overline\theta$ over $\mathbb{F}_p$. Then,
$$\deg(G)\le\deg(g)\le n-1<\deg(F).$$
Since $\overline{F(\theta)}=\overline0$, we have $F\equiv G\cdot\frac{F}{G}\pmod{\pi}$. Since $K$ is unramified, the uniformizer $\pi$ is $p$ multiplied by some unit in $K$, so we have $F\equiv G\cdot\frac{F}{G}\pmod{p}$ is reducible modulo $p$.
\end{proof}
\end{theorem}

\begin{remark}
We can see from the proof of Theorem \ref{th-2.3} that if we want to drop the unramifiedness assumption in this theorem, then the modulus will be $\pi$ instead of $p$. However, finding a uniformizer requires additional computation. Hence we do not use this stronger result.
\end{remark}

\begin{corollary}\label{co-2.6}
Let $K=\mathbb{Q}_p(\theta)$ be an unramified extension field of degree $n$ over $\mathbb{Q}_p$ with $\left|\theta\right|=1$. Let $G$ be the minimal polynomial of $\overline\theta$ over $\mathbb{F}_p$. If $\deg(G)=m$, then $1,\theta,\theta^2,\dots,\theta^{m-1}$ is an orthogonal basis of $V$ over $\mathbb{Q}_p$, where $V$ is the vector space generated by $1,\theta,\theta^2,\dots,\theta^{m-1}$ over $\mathbb{Q}_p$.
\begin{proof}
Similar to Theorem \ref{th-2.3}.
\end{proof}
\end{corollary}

We now begin to construct an orthogonal basis for a ramified extension field. The following lemma is an obvious corollary of Proposition 3.3 in \cite{ref-2}. It can be concluded by induction.

\begin{lemma}\label{le-2.7}
Let $K$ be an extension field of over $\mathbb{Q}_p$. Let $V_{i}\subset K$ be a vector space over $\mathbb{Q}_p$ of finite dimension $n_{i}>0$, $1\le i\le s$. Let ${\bm \alpha}_{i1},{\bm \alpha}_{i2},\dots,{\bm \alpha}_{in_{i}}$ be an orthogonal basis of $V_{i}$ over $\mathbb{Q}_p$. If
$$\left\{\left|{\bm v}_{i}\right|\Big|{\bm v}_{i}\in V_{i}\right\}\cap\left\{\left|{\bm v}_{j}\right|\Big|{\bm v}_{j}\in V_{j}\right\}=\{0\}$$
for all $1\le i<j\le s$. Then ${\bm \alpha}_{11},{\bm \alpha}_{12},\dots,{\bm \alpha}_{1n_{1}},\dots,{\bm \alpha}_{s1},{\bm \alpha}_{s2},\dots,{\bm \alpha}_{sn_{s}}$ is an orthogonal basis of $V=\bigoplus^{s}_{i=1}V_{i}$ over $\mathbb{Q}_p$.
\end{lemma}

\begin{theorem}\label{th-2.8}
Let $K$ be an extension field of degree $n$ over $\mathbb{Q}_p$. Let $f$ and $e$ be the residue degree and ramification index respectively. Let $\pi$ be a uniformizer of $K$ and $(s_{i})_{1\le i\le f}$ be a family in $R$ such that the image $\overline{s_{i}}\in k$ make up a basis of $k$ over $\mathbb{F}_p$. Then the family
$$(s_{i}\pi^{j})_{1\le i\le f,\ 0\le j\le e-1}$$
is an orthogonal basis of $K$ over $\mathbb{Q}_p$.
\begin{proof}
By Theorem \ref{th-0.4}, $n=ef$. We can prove that the elements in this family are linearly independent over $\mathbb{Q}_p$ (see  \cite{ref-4} page 99), so it is a basis of $K$ over $\mathbb{Q}_p$. Let $V_{j}$ be the vector space generated by $(s_i\pi^j)_{1\le i\le f}$ over $\mathbb{Q}_p$, $0\le j\le e-1$. Then $K=\bigoplus^{e-1}_{j=0}V_{j}$. Since $(\overline{s_{i}})_{1\le i\le f}$ are linearly independent over $\mathbb{F}_p$, by Theorem \ref{th-1.3}, $(s_i\pi^j)_{1\le i\le f}$ is an orthogonal basis of $V_{j}$ over $\mathbb{Q}_p$. Since $\left|\pi\right|=p^{-\frac{1}{e}}$, we have
$$\{\left|{\bm v}_{j}\right||{\bm v}_{j}\in V_{j}\}=\{0\}\cup p^{\mathbb{Z}-\frac{j}{e}}.$$
Then by Lemma \ref{le-2.7}, the family
$$(s_{i}\pi^{j})_{1\le i\le f,\ 0\le j\le e-1}$$
is an orthogonal basis of $K$ over $\mathbb{Q}_p$.
\end{proof}
\end{theorem}

\begin{example}
Let $K=\mathbb{Q}_3(\sqrt3+i)=\mathbb{Q}_3(\sqrt3,i)$ where $i^2=-1$. Then $K$ is an extension field of degree $n=4$ over $\mathbb{Q}_3$. Its residue degree is $f=2$ and ramification index is $e=2$. Since $\sqrt3$ is a uniformizer of $K$ and $1,i$ are linearly independent over $\mathbb{F}_3$, $\{1,i,\sqrt3,\sqrt3i\}$ is an orthogonal basis of $K$ over $\mathbb{Q}_3$.
\end{example}

\begin{lemma}\label{le-2.10}
Let $K_{1}=\mathbb{Q}_p(\theta)$ be an unramified extension field of degree $f$ over $\mathbb{Q}_p$ and $K_{2}=\mathbb{Q}_p(\pi)$ be a totally ramified extension field of degree $e$ over $\mathbb{Q}_p$. Then $K=\mathbb{Q}_p(\theta,\pi)$ is an extension field of degree $n=ef$ over $\mathbb{Q}_p$. Its residue degree is $f$ and ramification index is $e$.
\begin{proof}
Since its residue degree is at least $f$ and ramification index is at least $e$ and $[K:\mathbb{Q}_p]\le ef$, we have $[K:\mathbb{Q}_p]=ef$. Hence its residue degree is $f$ and ramification index is $e$.
\end{proof}
\end{lemma}

Finally, we can present our crucial result for constructing orthogonal basis.

\begin{theorem}\label{th-2.11}
Let $K_{1}=\mathbb{Q}_p(\theta)$ be an unramified extension field of degree $f$ over $\mathbb{Q}_p$ with $\left|\theta\right|=1$. Let $F$ be the minimal polynomial of $\theta$ over $\mathbb{Q}_p$. Let $K_{2}=\mathbb{Q}_p(\pi)$ be a totally ramified extension field of degree $e$ over $\mathbb{Q}_p$. Assume that $F$ is irreducible modulo $p$ and $\pi$ is a uniformizer of $K_{2}$. Then the family
$$(\theta^{i}\pi^{j})_{0\le i\le f-1,\ 0\le j\le e-1}$$
is an orthogonal basis of $K=\mathbb{Q}_p(\theta,\pi)$ over $\mathbb{Q}_p$.
\begin{proof}
By Lemma \ref{le-2.10}, $K$ is an extension field of degree $n=ef$ over $\mathbb{Q}_p$. Its residue degree is $f$ and ramification index is $e$. By Theorem \ref{th-2.3}, $1,\theta,\theta^{2},\dots,\theta^{f-1}$ is an orthogonal basis of $K_{1}$ over $\mathbb{Q}_p$. Then by Theorem \ref{th-1.3}, their image in $k$ are linearly indenpendent over $\mathbb{F}_p$ and therefore make up a basis of $k$ over $\mathbb{F}_p$. Finally, by Theorom \ref{th-2.8}, the family
$$(\theta^{i}\pi^{j})_{0\le i\le f-1,\ 0\le j\le e-1}$$
is an orthogonal basis of $K=\mathbb{Q}_p(\theta,\pi)$ over $\mathbb{Q}_p$.
\end{proof}
\end{theorem}

\section{Constructing Orthogonal Basis with Roots of Unit}\label{se-5}

Notice that roots of unity are helpful for us to construct an extension field which satisfies the assumption in Theorem \ref{th-2.3}.

\begin{theorem}\label{th-3.1}
Let $K=\mathbb{Q}_p(\theta)$ be an extension field of degree $n$ over $\mathbb{Q}_p$ and $\theta$ is a root of unity of order prime to $p$. Then $1,\theta,\theta^2,\dots,\theta^{n-1}$ is an orthogonal basis of $K$ over $\mathbb{Q}_p$.
\begin{proof}
Since a root of unity of order prime to $p$ generates over $\mathbb{Q}_p$ an unramified extension field (see \cite{ref-4} page 105), and the minimal polynomial of $\theta$ over $\mathbb{Q}_p$ is irreducible modulo $p$, the theorem is immediately from Theorem \ref{th-2.3}.
\end{proof}
\end{theorem}

\begin{example}
Let $K=\mathbb{Q}_3(i)$ where $i^2=-1$. Since $X^2+1$ is irreducible modulo $3$, $\{1,i\}$ is an orthogonal basis of $K$ over $\mathbb{Q}_3$.
\end{example}

In order to determine whether $F$ is irreducible modulo $p$, we need some results about cyclotomic cosets to factor $X^n-1$ modulo $p$ $(\gcd(n,p)=1)$.
\begin{definition}[cyclotomic coset]
Suppose $p$ is a prime number and $\gcd(n,p)=1$. For integer $s\in\{0,1,\dots,p-1\}$, the cyclotomic coset mod $n$ over $\mathbb{F}_p$ which contains integer $s$ is defined as
$$C_s=\{s,sp,sp^{2},\dots,sp^{m_s-1}\},$$
where $sp^{m_s}\equiv s\pmod{n}$, i.e., $m_s$ is the least positive integer satisfying $sp^{m_s}\equiv s\pmod{n}$.
\end{definition}

The following lemma is a well known result of cyclotomic polynomial. We omit the proof here.

\begin{lemma}[\cite{ref-17}]\label{le-3.2}
Let $p$ be a prime number. Let $\alpha$ be a primitive $n$th root of unity over $\mathbb{F}_p$. Let $s\in\{0,1,\dots,p-1\}$. Then
$$M^{(s)}(X)=\prod_{i\in C_s}(X-\alpha^i)\in\mathbb{F}_p[X]$$
is the minimal polynomial of $\alpha^s$ over $\mathbb{F}_p$. Moreover,
$$X^n-1=\prod_{s}M^{(s)}(X),$$
where $s$ runs over a set of cyclotomic coset representatives modulo $n$, gives the factorization of $X^n-1$ over $\mathbb{F}_p$.
\end{lemma}

\begin{lemma}\label{le-3.3}
Let $q$ be a prime number such that $q_{0}=\frac{q-1}{2}$ is also prime i.e., $q_0$ is a Sophie Germain prime. Let $p$ be an integer such that $p\not\equiv-1\pmod{q}$. Then $p$ is a primitive root modulo $q$ if and only if $p$ is not a quadratic residue modulo $q$.
\begin{proof}
$p$ is not a quadratic residue modulo $q$ if and only if $p^{\frac{q-1}{2}}\equiv-1\pmod{q}$. Since $p\not\equiv-1\pmod{q}$, this is equivalent to $p$ is a primitive root modulo $q$.
\end{proof}
\end{lemma}

Now we choose $q$, $q_{0}$ and $p$ as in Lemma \ref{le-3.3}. Assume $p$ is not a quadratic residue modulo $q$ so the order of $p$ modulo $q$ is $q-1$. Take $K_{1}=\mathbb{Q}_p(\theta)$ where $\theta$ is a primitive $q$th root of unity. By Lemma \ref{le-3.2}, the degree of the minimal polynomial of $\theta$ over $\mathbb{F}_p$ is $q-1$, hence it is also $q-1$ over $\mathbb{Q}_p$. Let $K_{2}=\mathbb{Q}_p(\pi)$ be a totally ramified extension field of degree $e$ over $\mathbb{Q}_p$.
Then by Theorem \ref{th-2.11}, the family
$$(\theta^{i}\pi^{j})_{0\le i\le q-2,\ 0\le j\le e-1}$$
is an orthogonal basis of $K=\mathbb{Q}_p(\theta,\pi)$ over $\mathbb{Q}_p$.\par
To make this orthogonal basis practical, we also need to find a primitive element of $K=\mathbb{Q}_p(\theta,\pi)$. It is well known that a finite separable extension is a simple extension, which is known as the primitive element theorem. If we consider this theorem in a field of characteristic $0$, then we only need the finiteness in the assumption. Moreover, there is a constructive proof of this theorem, see \cite{ref-3}(page 410). We summarize the proof as the following lemma.

\begin{lemma}\label{le-3.4}
Let $K=\mathbb{Q}_p(\theta,\pi)$ be a finite extension field over $\mathbb{Q}_p$. Let $f(X)$ and $g(X)$ be the minimal polynomial of $\theta$ and $\pi$ over $\mathbb{Q}_p$ respectively. Suppose their roots are $\theta^{(1)},\theta^{(2)},\dots,\theta^{(n)}$ and $\pi^{(1)},\pi^{(2)},\dots,\pi^{(m)}$ respectively. Choose $h\in\mathbb{Q}_p$ such that 
$$h\neq\frac{\pi^{(u)}-\pi^{(v)}}{\theta^{(s)}-\theta^{(t)}} \quad (1\le s,t\le n,\ 1\le u,v\le m).$$
Then $\zeta=h\theta+\pi$ is a primitive element of $K$ over $\mathbb{Q}_p$, i.e., $K=\mathbb{Q}_p(\zeta)$.
\end{lemma}

\begin{lemma}[\cite{ref-4} page 105]\label{le-3.5}
Let $K$ be any ultrametric extension of $\mathbb{Q}_p$. Then the distance between two distinct roots of unity of order prime to $p$ is $1$. 
\end{lemma}

\begin{lemma}\label{le-3.6}
Let $q\neq p$ be a prime number and $\theta$ be a primitive $q$th root of unity whose degree over $\mathbb{Q}_p$ is $f=q-1$. Let $G$ be an Eisenstein polynomial and $\pi$ be a root of $G(X)=0$. Let $K=\mathbb{Q}_p(\theta,\pi)$ and $\zeta=\theta+\pi$. Then $K=\mathbb{Q}_p(\zeta)$.
\begin{proof}
Let $\pi^{(1)},\pi^{(2)},\dots,\pi^{(e)}$ and $\theta^{(1)},\theta^{(2)},\dots,\theta^{(f)}$ be the conjugate elements of $\pi$ and $\theta$ over $\mathbb{Q}_p$ respectively, where $e=\deg{(F)}$. We have $\left|\pi^{(u)}\right|<1$. On the other hand, by Lemma \ref{le-3.5},
$$\left|\theta^{(s)}-\theta^{(t)}\right|=1 \quad (1\le s<t\le f).$$
Since
$$\left|\frac{\pi^{(u)}-\pi^{(v)}}{\theta^{(s)}-\theta^{(t)}}\right|=\left|\pi^{(u)}-\pi^{(v)}\right|<1,$$
we can take $h=1$ in Lemma \ref{le-3.4}. Therefore, $K=\mathbb{Q}_p(\theta+\pi)=\mathbb{Q}_p(\zeta)$.
\end{proof}
\end{lemma}

Finally, let $H$ be the minimal polynomial of $\theta$ over $\mathbb{Q}_p$. We can use resultant to obtain the minimal polynomial of $\zeta=\theta+\pi$, i.e., 
$$F(X)={\rm Res}_Y(G(Y),H(X-Y)).$$

We summarize the construction process as the following algorithm.\par\vskip 10pt

{\bf Algorithm} (constructing orthogonal basis with roots of unit).\par
{\bf Input:} two prime numbers $q$ and $q_{0}$ such that $q=2q_{0}+1$, a prime number $p$ such that $p\not\equiv-1\pmod{q}$ and $p$ is not a quadratic residue modulo $q$ and a positive integer $e$.\par
{\bf Output:} an extension field $K$ over $\mathbb{Q}_p$ of degree $n=(q-1)e$ and an orthogonal basis of $K$.
\begin{enumerate}
\item choose a random primitive $q$th root of unity $\theta$ and denote its minimal polynomial as $H$
\item choose a random Eisenstein polynomial $G$ of degree $e$ and choose a random root $\pi$ of $G(X)=0$
\item let $\zeta=\theta+\pi$, by Lemma \ref{le-3.6}, $\zeta$ is a primitive element of $\mathbb{Q}_p(\theta,\pi)$ over  $\mathbb{Q}_p$
\item let $F(X)={\rm Res}_Y(G(Y),H(X-Y))$
\item return $K=\mathbb{Q}_p(\zeta)$ (given by $F$) and an orthogonal basis
$$(\theta^{i}\pi^{j})_{0\le i\le q-2,\ 0\le j\le e-1},$$
where $\theta$ and $\pi$ are given by $H$ and $G$ respectively
\end{enumerate}\par

\section{Conclusion}

It is clear that the algorithm runs in polynomial time in $n$. The resultant in step 4 is an $(e+q-1)\times(e+q-1)$ determinant. If we take $e\approx q-1\approx\sqrt{n}$, then the storage requirement is only $O(n)$ and the time complexity is $O(n^{1.5})$ if we use decomposition methods to compute the resultant. In the worst case, where $\{e,q-1\}=\{1,n\}$, the storage requirement is $O(n^2)$ and the time complexity is $O(n^3)$. Therefore, our method is better than the Round 2 Algorithm and the Round 4 Algoritm in storage requirement. Moreover, according to \cite{ref-19}, the time complexity of the Round 2 Algorithm is more than $O(n^4)$. Our method is also better than the Round 2 Algorithm in time complexity.\par
The results in this paper are helpful to modify the $p$-adic signature scheme and public-key encryption cryptosystem. Simply using $\zeta=\theta+\pi$ to generate the extension field $K$ may be not secure. An adversary can guess the residue degree $f$ and subtract a primitive $f$th root of unit $\theta^{\prime}$ form $\zeta$. If it happends that $\theta^{\prime}=\theta$, then the adversary will obtain the uniformizer $\pi$ and break the scheme. Maybe we need to find a more complex primitive element $\zeta$ while do not increase time complexity much. More efforts are needed to carry out secure schemes. The method of constructing orthogonal basis in $p$-adic fields may have other applications. It is worth for further study and there is much work to do.

\section*{Acknowledgements}

This work was supported by National Natural Science Foundation of China(No. 12271517).

\end{document}